\documentclass[11pt]{article}

\usepackage{amssymb,amsthm}
\usepackage{caption,graphicx}
\usepackage[text={7.0in,9.5in},centering]{geometry}
\usepackage[comma,numbers,square,sort&compress]{natbib}

\setcaptionmargin{0.25in}

\newtheorem{Lemma}{Lemma}[section]
\newtheorem{Proposition}[Lemma]{Proposition}
\newtheorem{Theorem}{Theorem}

\makeatletter\@addtoreset{equation}{section}\makeatother

\def\cosh{\mathop\mathrm{cosh}\nolimits}
\def\sech{\mathop\mathrm{sech}\nolimits}

\def\errfn{\mathop\mathrm{errfn}\nolimits}

\newcommand{\R}{\mathbb{R}}             % reals
\newcommand{\rmO}{\mathrm{O}}           % Landau O
\newcommand{\rmd}{\,\mathrm{d}}         % derivatives
\newcommand{\rme}{\mathrm{e}}           % Euler constant

\begin{document}

\title{Toward nonlinear stability of sources via a modified Burgers equation}
\author{Margaret Beck \and Toan Nguyen \and Bj\"orn Sandstede \and Kevin Zumbrun}

\author{%
Margaret Beck\\
Department of Mathematics\\
Boston University\\
Boston, MA~02215, USA
\and
Toan Nguyen\\
Division of Applied Mathematics\\
Brown University\\
Providence, RI~02912, USA
\and
Bj\"orn Sandstede\\
Division of Applied Mathematics\\
Brown University\\
Providence, RI~02912, USA
\and
Kevin Zumbrun\\
Department of Mathematics\\
Indiana University\\
Bloomington, IN~47405, USA
}

\date{\today}
\maketitle

\begin{abstract}
Coherent structures are solutions to reaction-diffusion systems that are time-periodic in an appropriate moving frame and spatially asymptotic at $x=\pm\infty$ to spatially periodic travelling waves. This paper is concerned with sources which are coherent structures for which the group velocities in the far field point away from the core. Sources actively select wave numbers and therefore often organize the overall dynamics in a spatially extended system. Determining their nonlinear stability properties is challenging as localized perturbations may lead to a non-localized response even on the linear level due to the outward transport. Using a modified Burgers equation as a model problem that captures some of the essential features of coherent structures, we show how this phenomenon can be analysed and nonlinear stability be established in this simpler context.
\end{abstract}

%%%%%%%%%%%%%%%%%%%%%%%%%%%%%%%%%%%%%%%%%%%%%%%%%%%%%%%%%%%%%%%%%%%%%%%%%%%%

\section{Introduction}

In this paper, we analyse the long-time dynamics of solutions to the Burgers-type equation
\begin{equation}\label{E:toy}
\phi_t + c\tanh\left(\frac{cx}{2}\right) \phi_x = \phi_{xx} + \phi_x^2,
\qquad c>0
\end{equation}
with small localized initial data, where $x\in\R$, $t>0$, and $\phi(x,t)$ is a scalar function. The key feature of this equation as opposed to the usual Burgers equation is that the characteristic speeds are $c>0$ at spatial infinity and $-c<0$ at spatial minus infinity: hence, transport is always directed away from the shock interface at $x=0$ and not towards $x=0$ as would be the case for the Lax shocks of the standard Burgers equation.

\begin{figure}
\centering\includegraphics{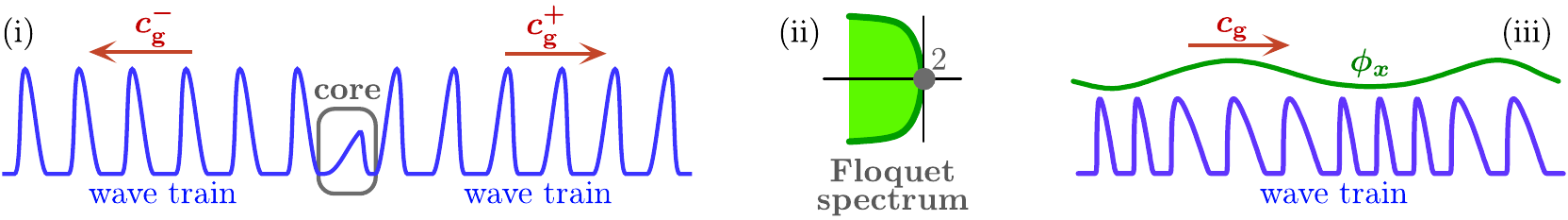}
\caption{Panel (i) shows the graph of a source $u^*(x,t)$ as a function of $x$ for fixed time $t$: the group velocities of the asymptotic wave trains point away from the core of the coherent structure. Panel~(ii) illustrates the Floquet spectrum of a spectrally stable source: the two eigenvalues at the origin correspond to translations in space and time, which, in contrast to the essential spectrum, cannot be moved by exponential weights. Panel~(iii) shows the behaviour of small phase $\phi$ or wave number $\phi_x$ perturbations of a wave train: to leading order, they are transported with speed given by the group velocity $c_\mathrm{g}$ without changing their shape \cite{DoelmanSSS}.}\label{f:1}
\end{figure}

We are interested in (\ref{E:toy}) due to its close connection with the dynamics of coherent structures that arise in reaction-diffusion systems
\begin{equation} \label{E:rd}
u_t = Du_{xx} + f(u), \qquad x\in\R, \qquad u\in\R^n.
\end{equation}
A coherent structure (or defect) is a solution $u^*(x,t)$ of (\ref{E:rd}) that is time-periodic in an appropriate moving frame $y=x-c^*t$ and spatially asymptotic to wave-train solutions, which are spatially periodic travelling waves of (\ref{E:rd}). Such structures have been observed in many experiments and in various reaction-diffusion models, and we refer to \cite{SandstedeScheel04_classify} for references and to Figure~\ref{f:1} for an illustration of typical defect profiles. For the sake of simplicity, we shall assume from now on that the speed $c^*$ of the defect we are interested in vanishes, so that the coherent structure is time-periodic. Coherent structures can be classified into several distinct types \cite{SaarloosHohenberg,Hecke,SandstedeScheel04_classify} that have different stability and multiplicity properties. This classification involves the group velocities of the asymptotic wave trains, and we therefore briefly review their definition and features. Wave trains of (\ref{E:rd}) are solutions of the form $u(x,t)=u_\mathrm{wt}(kx-\omega t;k)$, where the profile $u_\mathrm{wt}(y;k)$ is $2\pi$-periodic in the $y$-variable. Thus, $k$ and $\omega$ represent the spatial wave number and the temporal frequency, respectively, of the wave train. Wave trains typically exist as one-parameter families, where the frequency $\omega=\omega_\mathrm{nl}(k)$ is a function, the so-called nonlinear dispersion relation, of the wave number $k$, which varies in an open interval. The group velocity $c_\mathrm{g}$ of the wave train with wave number $k$ is defined as
\[
c_\mathrm{g} := \frac{\rmd \omega_\mathrm{nl}}{\rmd k}(k).
\]
The group velocity is important as it is the speed with which small localized perturbations of a wave train propagate as functions of time $t$, and we refer to Figure~\ref{f:1}(iii) for an illustration and to \cite{DoelmanSSS} for a rigorous justification of this statement. The classification of coherent structures mentioned above is based on the group velocities $c_\mathrm{g}^\pm$ of the asymptotic wave trains at $x=\pm\infty$. We are interested in sources for which $c_\mathrm{g}^-<0<c_\mathrm{g}^+$ as illustrated in Figure~\ref{f:1}(i) so that perturbations are transported away from the defect core towards infinity. Sources are important as they actively select wave numbers in oscillatory media; examples of sources are the Nozaki--Bekki holes of the complex Ginzburg--Landau equation.

From now on, we focus on a given source and discuss its stability properties with respect to the reaction-diffusion system (\ref{E:rd}). Spectral stability of a source can be investigated through the Floquet spectrum of the period map of the linearization of (\ref{E:rd}) about the time-periodic source. Spectral stability of sources was investigated in \cite{SandstedeScheel04_classify}, and we now summarize their findings. The Floquet spectrum of a spectrally stable source will look as indicated in Figure~\ref{f:1}(ii). A source $u^*(x,t)$ has two eigenvalues at the origin with eigenfunctions $u^*_x(x,t)$ and $u^*_t(x,t)$; the associated adjoint eigenfunctions are necessarily exponentially localized, so that the source has a well defined spatial position and temporal phase. There will also be two curves of essential spectrum that touch the origin and correspond to phase and wave number modulations of the two asymptotic wave trains. It turns out that the two eigenvalues at the origin cannot be removed by posing the linearized problem in exponentially weighted function spaces; the essential spectrum, on the other hand, can be moved to the left by allowing functions to grow exponentially at infinity.

\begin{figure}
\centering\includegraphics{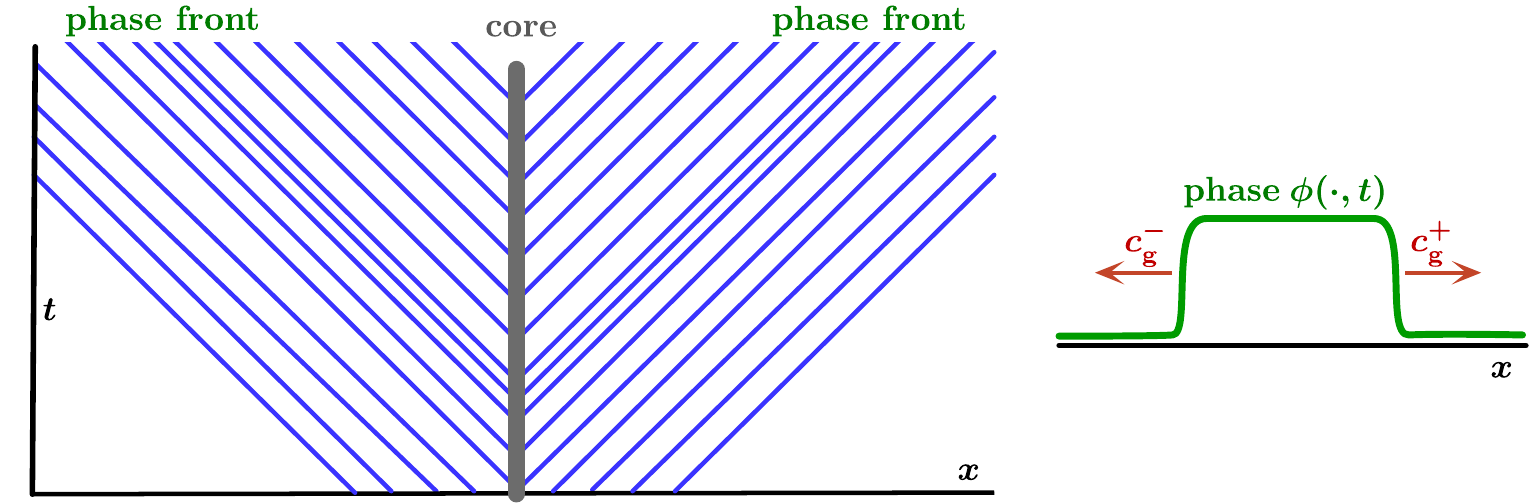}
\caption{The left panel contains a sketch of the space-time diagram of a perturbed source. The defect core will adjust in response to an imposed perturbation, and the emitted wave trains, whose maxima are indicated by the lines that emerge from the defect core, will therefore exhibit phase fronts that travel with the group velocities of the asymptotic wave trains away from the core towards $\pm\infty$. The right panel illustrates the profile of the anticipated phase function $\phi(x,t)$ defined in (\ref{e:phase}).}\label{f:2}
\end{figure}

The nonlinear stability of spectrally stable sources has not yet been established, and we now outline why this is a challenging problem. From a purely technical viewpoint, an obvious difficulty is related to the fact that there is no spectral gap between the essential spectrum and the imaginary axis. As discussed above, such a gap can be created by posing the linear problem on function spaces that contain exponentially growing functions, but the nonlinear terms will then not even be continuous. To see that these are not just technical obstacles, it is illuminating to discuss the anticipated dynamics near a source from an intuitive perspective. If a source is subjected to a localized perturbation, then one anticipated effect is that the defect core adjusts its position and its temporal phase in response. From its new position, the defect will continue to emit wave trains with the same selected wave number but there will now be a phase difference between the asymptotic wave trains at infinity and those newly emitted near the core. In other words, we expect to see two phase fronts that travel in opposite directions away from the core as illustrated in Figure~\ref{f:2}. The resulting phase dynamics can be captured by writing the perturbed solution $u(x,t)$ as
\begin{equation}\label{e:phase}
u(x,t) = u^*(x+\phi(x,t),t) + w(x,t),
\end{equation}
where we expect that the perturbation $w(x,t)$ of the defect profile decays in time, while the phase $\phi(x,t)$ resembles an expanding plateau as indicated in Figure~\ref{f:2} whose height depends on the initial perturbation through the spatio-temporal displacement of the defect core.

The preceding discussion indicates that localized perturbations of a defect will not stay localized but result instead in phase fronts that propagate towards infinity. As a first step towards a general nonlinear stability result for sources in reaction-diffusion systems, we focus here on the Burgers-type equation
\begin{equation}\label{e:tm}
\phi_t + c\tanh\left(\frac{cx}{2}\right) \phi_x = \phi_{xx} + \phi_x^2, \qquad c>0
\end{equation}
for the phase function $\phi(x,t)$ introduced in (\ref{e:phase}). The rationale for choosing this model problem is as follows. First, as established formally in \cite{HowardKopell} and proved rigorously in \cite{DoelmanSSS}, the integrated viscous Burgers equation captures the phase dynamics of wave trains over long time intervals. Since the Burgers equation does not admit sources, we add the artificial advection term on the left-hand side which creates the characteristic speeds $\pm c$ at $x=\pm\infty$ that account for the outgoing group velocities of the asymptotic wave trains of the source. While the inhomogeneous advection term can be thought of as fixing the position $x=0$ of the core, the equation still has a family of constant solutions, which correspond to different temporal phases of the underlying hypothetical sources. We therefore feel that gaining a detailed understanding of the long-time dynamics of (\ref{e:tm}) for small localized initial data will shed significant light on the expected dynamics of sources and on the techniques needed to analyse their stability. We emphasize that the dynamics of wave trains of reaction-diffusion systems under non-localized phase perturbations was investigated only recently in \cite{SSSUecker}; the methods used there rely on renormalization-group techniques which seem difficult to generalize to sources. On the other hand, the analysis presented here is currently far more limited in terms of the equations it applies to.

\begin{figure}
\centering\includegraphics{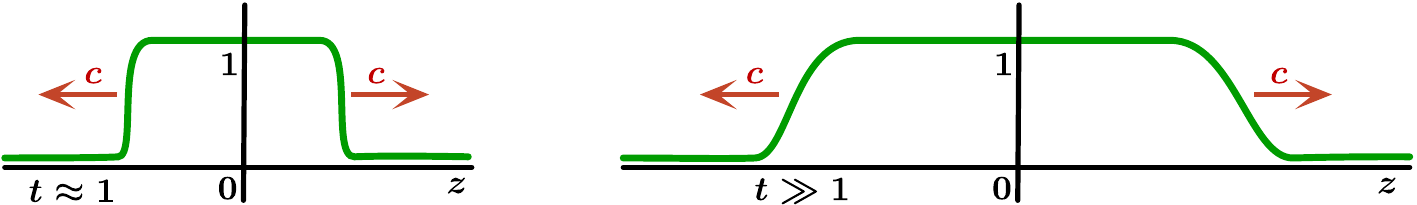}
\caption{Shown are the graphs of the function $\errfn((-z+ct)/\sqrt{4t})-\errfn((-z-ct)/\sqrt{4t})$ for smaller and larger values of $t$, which resemble plateaus of height approximately equal to one that spread outwards with speed $\pm c$.}\label{f:3}
\end{figure}

The linearization of (\ref{e:tm}) about $\phi=0$ is given by
\begin{equation}\label{E:toy_lin}
\phi_t = \phi_{xx} - c\tanh\left(\frac{cx}{2}\right) \phi_x.
\end{equation}
The spectrum of the operator in (\ref{E:toy_lin}) is as shown in Figure~\ref{f:1}(ii) except that there is only one embedded eigenvalue at the origin that cannot be moved by exponential weights: the eigenfunction of the eigenvalue at $\lambda=0$ is given by the constant function $\frac{c}{4}$, and the associated adjoint eigenfunction $\psi$ is given by
\[
\psi(y) = \sech^2\left(\frac{cy}{2}\right).
\]
Equation (\ref{E:toy_lin}) admits\footnote{Note that (\ref{E:toy_lin}) is the formal adjoint of the linearization of the standard viscous Burgers equation $u_t=u_{zz}-2uu_z$ about the Lax shock $\bar{u}(x)=(c/2)[1-\tanh(cx/2)]$ with $x=z-ct$, whose Green's function can be found via the linearized Cole--Hopf transformation by setting $w(x,t)=\cosh(cx/2)\int_{-\infty}^x u(y,t)\rmd y$. The Green's function of (\ref{E:toy_lin}) can then be constructed by reversing the roles of $x$ and $y$ in the Green's function for the Lax shock linearization.} the Green's function
\begin{eqnarray} \label{E:toy_greens}
\mathcal{G}(x,y,t) & = &
\frac{1}{\sqrt{4\pi t}} \rme^{-\frac{(x-y+ct)^2}{4t}}\frac{1}{1+\rme^{cy}}
+ \frac{1}{\sqrt{4\pi t}} \rme^{-\frac{(x-y-ct)^2}{4t}}\frac{1}{1+\rme^{-cy}}
\\ \nonumber 
& & \qquad + \frac{c}{4}\left[ \errfn\left(\frac{y-x+ct}{\sqrt{4t}}\right) - \errfn\left(\frac{y-x-ct}{\sqrt{4t}}\right) \right] \psi(y),
\end{eqnarray}
where the error function is given by
\[
\errfn(z) = \frac{1}{\sqrt{\pi}}\int_{-\infty}^z \rme^{-s^2} \rmd s.
\]
The first two terms in the Green's function are Gaussians that move with speed $\pm c$ away from the core and decay at the rate $1/\sqrt{t}$, while the term comprised of the difference of the two error functions produces a plateau of constant height $\frac{c}{4}$ that spreads outward as indicated in Figure~\ref{f:3}.

For sufficiently localized initial data $\phi_0(x)$, the solution $\phi(x,t)$ to the linear equation (\ref{E:toy_lin}) is therefore given by
\[
\phi(x,t) = \int_\mathbb{R} \mathcal{G}(x,y,t) \phi_0(y) \rmd y,
\]
which converges to a constant:
\[
\phi(x,t) \longrightarrow \frac{c}{4}\int_\mathbb{R}\psi(y) \phi_0(y) \rmd y \quad\mbox{as}\quad t\longrightarrow\infty
\]
for each fixed $x$. In particular, one cannot expect the solution associated with localized initial data to remain localized for all time. In fact, the same is true for the nonlinear equation (\ref{e:tm}): the Cole--Hopf transformation
\[
\tilde\phi(x,t) = \rme^{\phi(x,t)}-1, \qquad \phi(x,t) = \log\Big[1+\tilde\phi(x,t)\Big]
\]
relates solutions $\phi(x,t)$ of (\ref{e:tm}) and solutions $\tilde\phi(x,t)$ of the linearization (\ref{E:toy_lin}), and we conclude that solutions $\phi(x,t)$ of (\ref{e:tm}) are given by
\[
\phi(x,t) = \log\left[1+\int_\mathbb{R}\mathcal{G}(x,y,t)\phi_0(y)\rmd y \right].
\]
For $t\to\infty$, these solutions converge again pointwise in $x$ to the constant
\[
\log\left[ \int_\mathbb{R}\psi(y)\left(\rme^{\phi_0(y)} - 1\right)\rmd y + 1 \right].
\]
The Cole--Hopf transformation does not, however, extend to more general equations, and we therefore pursue here a different approach that, we hope, will also be useful when investigating the nonlinear stability of sources in general reaction-diffusion systems. 

In order to analyze the dynamics of (\ref{e:tm}) in a way that may also be applicable in the case of a general reaction-diffusion equation, we need to find an appropriate ansatz for small-amplitude solutions. To do so, we define $\mathcal{B}(x,t)$ by 
\[
\mathcal{B}(x,t) := \mathcal{G}(x,0,t+1),
\]
where $\mathcal{G}(x,y,t)$ is the Green's function defined in (\ref{E:toy_greens}), and note that
\begin{equation}\label{def-phia}
\phi^*(x,t,p) := \log\left( 1 + p \mathcal{B}(x,t) \right)
\end{equation}
is then a solution\footnote{In fact, $\phi^*(x,t,p)$ does not need to satisfy (\ref{e:tm}) exactly: our analysis goes through provided $\mathcal{B}(x,t)$ can be chosen such that $\phi^*(x,t,p)$ satisfies (\ref{e:tm}) approximately with an error of order of $(1+t)^{-1}$ times a Gaussian.} of (\ref{e:tm}) for each fixed $p\in\R$. We now investigate the long-time dynamics of solutions to (\ref{e:tm}) with initial data $\phi(x,0)=\phi_0(x)$ using the ansatz
\begin{equation}\label{E:ansatz}
\phi(x,t) = \phi^*(x,t,p(t)) + v(x,t),
\end{equation}
where $p(t)$ is a real-valued function, and $v(x,t)$ is a remainder term. At time $t=0$, we normalize the decomposition in (\ref{E:ansatz}) by choosing $p(0)=p_0$ such that
\[
\int_\mathbb{R} \psi(x) \big[\phi_0(x) - \phi^*(x,0,p_0)\big] \rmd x = 0.
\]
We will prove in \S\ref{S:proof} that a unique $p_0$ with this property exists for each sufficiently small localized initial condition $\phi_0$. Note that the leading behavior $  \phi^*(x,t,p(t))$ corresponds to a plateau with height $\log(1+p(t))$ of length $2ct$ that spreads outward with speed $\pm c$. The main result of this paper is as follows.

\begin{Theorem}\label{Thm:main}
For each $\gamma\in(0,\frac12)$, there exist constants $\epsilon_0,\eta_0,C_0,M_0>0$ such that the following is true. If $\phi_0\in C^1$ satisfies
\begin{equation}\label{e:phi0}
\epsilon := \|\rme^{x^2/M_0}\phi_0\|_{C^1} \leq \epsilon_0
\end{equation}
then the solution $\phi(x,t)$ of (\ref{e:tm}) with $\phi(\cdot,0)=\phi_0$ exists globally in time and can be written in the form
\[
\phi(x,t) =   \phi^*(x,t,p(t)) + v(x,t)
\]
for appropriate functions $p(t)$ and $v(x,t)$ with $  \phi^*$ as in (\ref{def-phia}). Furthermore, there is a $p_\infty\in\R$ with $|p_\infty|\leq C_0$ such that
\[
|p(t)-p_\infty| \leq \epsilon C_0 \rme^{-\eta_0 t},
\]
and $v(x,t)$ satisfies
\[
|v(x,t)| \leq \frac{\epsilon C_0}{(1+t)^{\gamma}}
\left( \rme^{-\frac{(x+ct)^2}{M_0(t+1)}} + \rme^{-\frac{(x-ct)^2}{M_0(t+1)}} \right),
\quad
|v_x(x,t)| \leq \frac{\epsilon C_0}{(1+t)^{\gamma+1/2}}
\left( \rme^{-\frac{(x+ct)^2}{M_0(t+1)}} + \rme^{-\frac{(x-ct)^2}{M_0(t+1)}} \right)
\]
for all $t\geq0$. In particular, $\|v(\cdot,t)\|_{L^r}\to0$ as $t\to\infty$ for each fixed $r>\frac{1}{2\gamma}$.
\end{Theorem}

Note that our result assumes that the initial condition is strongly localized in space. We believe that this assumption can be relaxed significantly; for the purposes of this paper, however, phase fronts are created even by highly localized initial data, and the key difficulties are therefore present already in the more specialized situation of Theorem~\ref{Thm:main}.

The exponential convergence of $p(t)$ reflects the intuition that a source has a well-defined position due to the exponential localization of the adjoint eigenfunction $\psi$. The asymptotics of the perturbation $v$ is given by moving Gaussians that decay only like $(1+t)^{-\gamma}$ for each fixed $\gamma\in(0,\frac12)$, rather than with $(1+t)^{-\frac12}$ as expected from the dynamics of the viscous Burgers equation. This weaker result is due to the form (\ref{E:ansatz}) of our ansatz, which effectively creates a nonlinear term that is proportional to $g(x)uu_x$ for some function $g(x)$. While this term resembles the term $2uu_x$ in Burgers equation, it does not respect the conservation-law structure. Thus, we only obtain decay at the above rate. Although this may not be optimal, it allows us to avoid terms that grow logarithmically in the nonlinear iteration in \S\ref{S:proof}. We do not know if it is possible to improve this rate to $\gamma=\frac12$ by adjusting our ansatz appropriately.

The remainder of this paper, \S\ref{S:proof}, is devoted to the proof of Theorem~\ref{Thm:main}.

%%%%%%%%%%%%%%%%%%%%%%%%%%%%%%%%%%%%%%%%%%%%%%%%%%%%%%%%%%%%%%%%%%%%%%%%%%%%

\section{Proof of the main theorem} \label{S:proof}

To prove Theorem~\ref{Thm:main}, we set up integral equations for $(p,v)$ and solve them using a nonlinear iteration scheme. Throughout the proof, we denote by $C$ possibly different positive constants that depend only on the underlying equation but not on the initial data or on space or time.

\subsection{Derivation of an integral formulation}

Substituting the ansatz
\[
\phi(x,t) = \log\left( 1+p(t)\mathcal{B}(x,t) \right) + v(x,t), \qquad
\mathcal{B}(x,t) := \mathcal{G}(x,0,t+1),
\]
into the Burgers-type equation
\[
\phi_t + c\tanh\left(\frac{cx}{2}\right) \phi_x = \phi_{xx} + \phi_x^2,
\]
we find that $(p,v)$ needs to satisfy the equation
\begin{equation}\label{E:v}
v_t = v_{xx} - c\tanh\left(\frac{cx}{2}\right) v_x + v_x^2 - \frac{\dot{p}}{1+\frac{c }{4}p} \mathcal{G}(x,0,t+1)+ \mathcal{N}(x,t,p,\dot{p},v_x),
\end{equation}
where the nonlinear function $\mathcal{N}$ is given by
\[
\mathcal{N}(x,t,p,\dot{p},v_x) := \frac{2pv_x\mathcal{B}_x(x,t)}{1+p\mathcal{B}(x,t)} + \dot{p} \left( \frac{\mathcal{B}(x,t)}{1+\frac{c}{4}p} - \frac{\mathcal{B}(x,t)}{1+p\mathcal{B}(x,t)} \right).
\]
The idea of the proof is to use an appropriate integral representation for $v$ that will allow us to set up a nonlinear iteration argument to show that solutions $(p,v)$ of (\ref{E:v}) exist and that they satisfy the desired decay estimates in space and time. Recall from (\ref{E:toy_greens}) the expression
\begin{eqnarray} \label{e:gf}
\mathcal{G}(x,y,t) & = &
\frac{1}{\sqrt{4\pi t}} \rme^{-\frac{(x-y+ct)^2}{4t}}\frac{1}{1+\rme^{cy}}
+ \frac{1}{\sqrt{4\pi t}} \rme^{-\frac{(x-y-ct)^2}{4t}}\frac{1}{1+\rme^{-cy}}
\\ \nonumber 
& & + \frac{c}{4}\left( \errfn\left(\frac{y-x+ct}{\sqrt{4t}}\right) - \errfn\left(\frac{y-x-ct}{\sqrt{4t}}\right) \right) \psi(y)
\end{eqnarray}
for the Green's function of the linear problem (\ref{E:toy_lin}) and note that it satisfies the identity
\[
\int_\mathbb{R} \mathcal{G}(x,y,t-s)\mathcal{G}(y,0,s+1)\rmd y = \mathcal{G} (x,0,t+1).
\]
Using this identity and the variation-of-constants formula, we can rewrite (\ref{E:v}) in integral form and obtain\footnote{We shall often use the notation $[f+g](x,t)$ to denote $f(x,t)+g(x,t)$.}
\begin{eqnarray}\label{E:v_init}
v(x,t) & = &- \frac{4}{c} \mathcal{G}(x,0,t+1)\left( \log\left(1+\frac{c p(t)}{4} \right) - \log\left(1+\frac{c p(0)}{4}\right) \right) \\ \nonumber &&
+  \int_\mathbb{R}\mathcal{G}(x,y,t)v_0(y)\rmd y  
+  \int_0^t \int_\mathbb{R} \mathcal{G}(x,y,t-s) \left[ v_y^2 +\mathcal{N}(\cdot,\cdot,p,\dot{p},v_y)\right](y,s) \rmd y \rmd s.
\end{eqnarray}
The expression (\ref{e:gf}) of the Green's function $\mathcal{G}(x,y,t)$ shows that the terms involving the error functions do not provide any temporal decay. We need to treat these terms separately and therefore write the Green's function $\mathcal{G}(x,y,t)$ as
\[
\mathcal{G}(x,y,t) = \mathcal{E}(x,y,t) + \tilde{\mathcal{G}}(x,y,t),
\]
where
\[
\mathcal{E}(x,y,t) = e(x,t)\psi(y), \qquad
e(x,t) := \frac{c}{4}\left( \errfn\left(\frac{x+ct}{\sqrt{4t}}\right) - \errfn\left(\frac{x-ct}{\sqrt{4t}}\right) \right)
\]
and 
\begin{eqnarray}\label{def-tG}
\tilde{\mathcal{G}}(x,y,t) & := & \mathcal{G}(x,y,t) - \mathcal{E}(x,y,t)
\\ & = & \nonumber
\frac{1}{\sqrt{4\pi t}} \rme^{-\frac{(x-y+ct)^2}{4t}}\frac{1}{1+\rme^{cy}}
+ \frac{1}{\sqrt{4\pi t}} \rme^{-\frac{(x-y-ct)^2}{4t}}\frac{1}{1+\rme^{-cy}}
\\ && \nonumber
+ \left( \errfn\left(\frac{y-x+ct}{\sqrt{4t}}\right) - \errfn\left(\frac{y-x-ct}{\sqrt{4t}}\right) \right) \frac{c}{4}\sech^2(\frac{cy}{2})
\\ && \nonumber
- \left( \errfn\left(\frac{-x+ct}{\sqrt{4t}}\right) - \errfn\left(\frac{-x-ct}{\sqrt{4t}}\right) \right) \frac{c}{4}\sech^2(\frac{cy}{2}).
\end{eqnarray}
A calculation shows that
\[
\left| \errfn\left(\frac{y-x\pm ct}{\sqrt{4t}}\right) - \errfn\left(\frac{-x\pm ct}{\sqrt{4t}}\right) \right| \frac{c}{4} \sech^2(\frac{cy}{2}) \leq
C t^{-1/2} \left(\rme^{-\frac{(x-y+ct)^2}{4t}} + \rme^{-\frac{(x-y-ct)^2}{4t}} \right) \rme^{-c|y|/4},
\]
and we conclude that
\begin{equation}\label{e:tildeg}
|\tilde{\mathcal{G}}(x,y,t)| \leq C t^{-1/2} \left( \rme^{-\frac{(x-y+ct)^2}{4t}} + \rme^{-\frac{(x-y-ct)^2}{4t}} \right).
\end{equation}
To analyse (\ref{E:v_init}), we consider the initial condition $\phi_0(x)=\phi^*(x,0,p(0))+v_0(x)$ and show that, if $\phi_0 $ is sufficiently small in $L^\infty$, then $p(0)=p_0$ can be chosen such that
\[
\int_\mathbb{R} \psi(y) v_0(y) \rmd y = 0
\]
or, equivalently,
\begin{equation}\label{choice-p0}
\int_\mathbb{R} \psi(y) [\phi_0(y) -   \phi^*(y,0,p_0)] \rmd y = 0.
\end{equation}
To prove this claim, we observe that
\begin{equation}\label{choice-p1}
\phi^*(y,0,p_0) = \log(1 + p_0 \mathcal{B}(y,0)) = p_0 \mathcal{G}(y,0,1) + \rmO(p_0^2).
\end{equation} We note that the term $\rmO(p_0^2)$ in (\ref{choice-p1}) is bounded uniformly in $y\in\R$, and substitution into (\ref{choice-p0}) therefore gives the equation
\[
\int_\mathbb{R} \psi(y) \phi_0(y) \rmd y = p_0 \int_\mathbb{R} \psi(y) \mathcal{G}(y,0,1)\rmd y + \rmO(p_0^2),
\]
which can be solved uniquely for $p_0=p(0)$ near zero for each $\phi_0\in L^\infty$ for which $\| \phi_0\|_{L^\infty}$ is small enough. In particular, there is a constant $C>0$ such that the resulting initial value $p(0)$ satisfies 
\[
|p(0)| \leq C \| \phi_0\|_{L^\infty} \leq C\epsilon.
\]
We will construct $p(t)$ such that
\begin{equation}\label{eqs-dotp}
\dot{p}(t) = \left(1+ \frac{cp(t)}{4} \right) \int_\mathbb{R} \psi(y)\left[ v_y^2 +\mathcal{N}(\cdot,\cdot,p,\dot{p},v_y)\right](y,t) \rmd y
\end{equation}
or, equivalently, that
\begin{equation}\label{eqs-p}
\log\left(1+\frac{c p(t)}{4}\right) = \log\left(1+\frac{cp_0}{4}\right) +\frac{c}{4} \int_0^t \int_\mathbb{R} \psi(y)\left[ v_y^2 +\mathcal{N}(\cdot,\cdot,p,\dot{p},v_y)\right](y,s) \rmd y \rmd s
\end{equation}
for all $t\geq0$. Substituting (\ref{eqs-p}) into (\ref{E:v_init}), we obtain the equation
\begin{eqnarray}\label{E:v_p}
v(x,t) & = & -\frac{4}{c}\tilde{\mathcal{G}}(x,0,t+1) \left( \log\left(1+\frac{c p(t)}{4}\right) - \log\left(1+\frac{c p_0}{4}\right) \right) + \int_\mathbb{R}\tilde{\mathcal{G}}(x,y,t)v_0(y)\rmd y  \\ \nonumber &&
+ \int_0^t \int_\mathbb{R} \tilde{\mathcal{G}}(x,y,t-s) \left[ v_y^2 +\mathcal{N}(\cdot,\cdot,p,\dot{p},v_y) \right](y,s) \rmd y \rmd s
\\ \nonumber &&
+  \int_0^t \int_\mathbb{R} \left(e(x,t-s) - e(x,t+1)\right) \psi(y) \left[ v_y^2 +\mathcal{N}(\cdot,\cdot,p,\dot{p},v_y)\right](y,s) \rmd y \rmd s
\end{eqnarray}
for $v(x,t)$.

To analyse solutions of (\ref{eqs-dotp}) and (\ref{E:v_p}), we introduce template functions that capture the anticipated spatio-temporal behaviour of solutions. For each fixed choice of $\gamma\in(0,\frac12)$ and $M>0$, we define
\begin{equation} \label{E:def_theta}
\theta_1(x,t) = \frac{1}{(1+t)^{\gamma}} \left( \rme^{-\frac{(x-ct)^2}{M(t+1)}} + \rme^{-\frac{(x+ct)^2}{M(t+1)}}\right), \qquad
\theta_2(x,t) = \frac{1}{(1+t)^{\gamma+1/2}} \left( \rme^{-\frac{(x-ct)^2}{M(t+1)}} + \rme^{-\frac{(x+ct)^2}{M(t+1)}}\right)
\end{equation}
and let
\[
h_1(t) := \sup_{y\in\mathbb{R}, 0\leq s\leq t} \left[\frac{|v|}{\theta_1}+\frac{|v_y|}{\theta_2} \right](y,s), \qquad
h_2(t) := \sup_{0\leq s\leq t} |\dot{p}(s)|\rme^{c^2 s/M}, \qquad
h(t) := h_1(t) + h_2(t).
\]
We will later choose $M\gg1$.

We remark that we do know existence and smoothness of $(v,p)$ for short times: Indeed, we can solve the original PDE for $\phi(x,t)$ for short times and can substitute the resulting expression into (\ref{eqs-dotp}) upon using (\ref{E:ansatz}). The resulting integral equation has a solution $\dot{p}(t)$ for small times and, using again (\ref{E:ansatz}), we find a smooth function $v$ that then satisfies (\ref{E:v_p}). Furthermore, using that $\phi_0$ satisfies (\ref{e:phi0}) by assumption, we see that $h(t)$ is well defined and continuous for $0<t\ll1$. Finally, standard parabolic theory implies that $h(t)$ retains these properties as long as $h(t)$ stays bounded. The key issue is therefore to show that $h(t)$ stays bounded for all times $t>0$, and this is what the following proposition asserts.

\begin{Proposition}\label{P:claim-zeta}
For each $\gamma\in(0,\frac12)$, there exist positive constants $\epsilon_0,C_0,M$ such that
\begin{equation}\label{claim-zeta}
h_1(t) \leq C_0 (\epsilon  + h_2(t) + h(t)^2), \qquad
h_2(t) \leq C_0 (\epsilon + h(t)^2).
\end{equation}
for all $t\geq0$ and all initial data $u_0$ with $\epsilon:=\|\rme^{x^2/M}\phi_0\|_{C^1}\leq\epsilon_0$.
\end{Proposition}

Using this proposition, we can add the inequalities in (\ref{claim-zeta}) and eliminate $h_2$ on the right-hand side to obtain
\[
h(t) \leq C_0(C_0+1)(\epsilon + h(t)^2).
\]
Using this inequality and the continuity of $h(t)$, we find that $h(t)\leq 2C_0(C_0+1)\epsilon $ for all $t\geq0$ provided $0<\epsilon\leq\epsilon_0$ is sufficiently small. Thus, Theorem~\ref{Thm:main} will be proved once we establish Proposition~\ref{P:claim-zeta}. The following sections will be devoted to proving this proposition.

%%%%%%%%%%%%%%%%%%%%%%%%%%%%%%%%%%%%%%%%%%%%%%%%

\subsection{Estimates of the nonlinear term}

We begin by deriving estimates of the nonlinear term
\[
\mathcal{N}(x,t,p,\dot{p},v_x) =  \frac{2pv_x\mathcal{B}_x(x,t)}{1+p\mathcal{B}(x,t)} + \dot{p} \left( \frac{\mathcal{B}(x,t)}{1+\frac{c}{4}p} - \frac{\mathcal{B}(x,t)}{1+p\mathcal{B}(x,t)} \right)
\]
that appears in (\ref{eqs-dotp}) and (\ref{E:v_p}). Recalling that
\[
e(x,t) = \frac{c}{4}\left( \errfn\left(\frac{x+ct}{\sqrt{4t}}\right) - \errfn\left(\frac{x-ct}{\sqrt{4t}}\right) \right)
\]
and using the definition of the error function, we see that
\[
\left| e(x,t) \left(\frac{c}{4} - e(x,t)\right) \right| \leq C \left(\rme^{-\frac{(x + ct)^2}{8t}} +\rme^{-\frac{(x - ct)^2}{8t}} \right).
\]
Since $\mathcal{B}(x,t)=e(x,t+1)+\tilde{\mathcal{G}}(x,0,t+1)$, we conclude from (\ref{e:tildeg}) that
\[
\left|\mathcal{B}(x,t) \left(\frac{c}{4}-\mathcal{B}(x,t)\right) \right| \leq C \left(\rme^{-\frac{(x + ct)^2}{8(t+1)}} +\rme^{-\frac{(x - ct)^2}{8(t+1)}} \right)
\]
and thus
\[
\frac{\mathcal{B}(x,t)}{1+p\mathcal{B}(x,t)} - \frac{\mathcal{B}(x,t)}{1+\frac{c}{4}p} \leq C|p| \left(\rme^{-\frac{(x + ct)^2}{8(t+1)}} +\rme^{-\frac{(x - ct)^2}{8(t+1)}} \right).
\]
In addition, we have
\[
\left| \frac{\mathcal{B}_x(x,t)}{1+p\mathcal{B}(x,t)} \right| \leq C(1+t)^{-1/2} \left(\rme^{-\frac{(x + ct)^2}{8(t+1)}} + \rme^{-\frac{(x - ct)^2}{8(t+1)}} \right).
\]
Combining these estimates, we therefore obtain
\begin{equation}\label{bound-N}
|\mathcal{N}(x,t,p,\dot{p},v_x)| \leq C \left( (1+t)^{-1/2} |p| |v_x| + |p\dot{p}| \right) \left(\rme^{-\frac{(x + ct)^2}{8(t+1)}} +\rme^{-\frac{(x - ct)^2}{8(t+1)}} \right)
\end{equation}
uniformly in $x\in\mathbb{R}$ and $t\geq0$, provided $p$ is sufficiently small.

\subsection{Estimates for $h_2(t)$}

To establish the claimed estimate for $h_2(t)$, recall from (\ref{eqs-dotp}) that 
\begin{equation}\label{est-qdot1}
|\dot{p}(t)| \leq  C (1+|p(t)|) \int_\mathbb{R} \psi(y) \Big| \left[ v_y^2 +  \mathcal{N}(\cdot,\cdot,p,\dot{p},v_y) \right](y,t) \Big| \rmd y. 
\end{equation}
We shall show that there is a constant $C_1=C_1(M)$ such that
\begin{equation}\label{c:1}
|\dot{p}(t)|\leq C_1 \rme^{-c^2 t/M} (\epsilon + h^2(t)),
\end{equation}
which then establishes the estimate for $h_2(t)$ stated in Proposition~\ref{P:claim-zeta}. To show (\ref{c:1}), we note that the definitions of $p(t)$ and $h_2(t)$ imply that
\begin{equation}\label{p-estimate}
|p(t)| \leq |p(0)| + \int_0^t |\dot{p}(s)| \rmd s
\leq |p(0)| + \int_0^t \rme^{-c^2 s/M} h_2(t) \rmd s
\leq C_1 (\epsilon + h_2(t))
\end{equation}
and therefore
\begin{equation} \label{E:p_est}
|p(t)\dot{p}(t)| \leq C_1 \rme^{-c^2 t/M} h_2(t)(\epsilon + h_2(t)) \leq C_1 \rme^{-c^2 t/M}(\epsilon + h(t)^2).
\end{equation}  
Next, we use the estimate $|\psi(y)|\leq 2\rme^{-c|y|}$, the bound (\ref{bound-N}) on $\mathcal{N}$, and the inequality
\[
\rme^{-\frac{c|y|}{2}} \rme^{-\frac{(y\pm ct)^2}{M(1+t)}} \leq C_1 \rme^{-\frac{c|y|}{4}} \rme^{-\frac{c^2 t}{M}},
\]
which holds for each $M\geq8$, to conclude that
\begin{equation}\label{est-phiN}
|\psi(y)v_y^2(y,t)| \leq \frac{C_1}{(1+t)^{1+2\gamma}} C \rme^{-\frac{c}{2} |y|} \rme^{-\frac{2c^2}{M}t}  h_1(t)^2
\end{equation}
and therefore
\begin{eqnarray}
|\psi(y)\mathcal{N}(y,t,p,\dot{p},v_y)| & \leq & C_1 (1+t)^{-1/2} \left( \rme^{-\frac{(x + ct)^2}{8(t+1)}} + \rme^{-\frac{(x - ct)^2}{8(t+1)}} \right) \psi (y) |v_y(y,t)| |p(t)|
\nonumber \\ \nonumber & &
+ C_1  \left( \rme^{-\frac{(x + ct)^2}{8(t+1)}} + \rme^{-\frac{(x - ct)^2}{8(t+1)}} \right) \psi(y)|p(t)\dot{p}(t)|
\\ \nonumber & \leq &
C_1 (1+t)^{-1/2} \rme^{-\frac{c}{2}|y| - \frac{2c^2}{M} t} h_1(t) (\epsilon + h_2(t)) + C_1 \rme^{-\frac{c}{2}|y| - \frac{2c^2}{M} t}(\epsilon + h(t)^2) \\ \label{est-phiN3} & \leq &
C_1\rme^{-\frac{c}{2}|y| - \frac{2c^2}{M} t} (\epsilon + h(t)^2).
\end{eqnarray}
Using these estimates in (\ref{est-qdot1}), we arrive readily at (\ref{c:1}), thus proving the estimate for $h_2(t)$ stated in Proposition~\ref{P:claim-zeta}.

%%%%%%%%%%%%%%%%%%%%%%%%%%%%%%%%%%%%%%%%%%%%%%%%%%

\subsection{Estimates for $v(x,t)$ and $v_x(x,t)$}\label{sec-estv}

In this section, we will establish the pointwise bounds
\begin{eqnarray}
|v(x,t)|   & \leq & C_1 (\epsilon + h_2(t) + h(t)^2) \theta_1(x,t)
\label{key-est-phi} \\ \label{key-est-phix}
|v_x(x,t)| & \leq & C_1 (\epsilon + h_2(t) + h(t)^2) \theta_2(x,t)
\end{eqnarray}
for $v(x,t)$ and $v_x(x,t)$, respectively, which taken together prove the inequality for $h_1(t)$ stated in Proposition~\ref{P:claim-zeta}. In particular, the proof of Proposition~\ref{P:claim-zeta} is complete once the two estimates above are established. We denote by $C_1$ possibly different constants that depend only on the choice of $\gamma,M$ so that $C_1=C_1(\gamma,M)$.

We focus first on $v(x,t)$. From the integral formulation (\ref{E:v_p}) for $v(x,t)$, we find that
\begin{eqnarray}\label{phi-estimate}
|v(x,t)| & \leq &
C_1\left(|p_0| + |p(t)|\right)\tilde{\mathcal{G}}(x,0,t+1) + \int_\mathbb{R} \tilde{\mathcal{G}}(x,y,t)|v_0(y)|\rmd y 
\\ \nonumber & &
+ \int_0^t \int_\mathbb{R} \Big| \tilde{\mathcal{G}}(x,y,t-s) \left[ v_y^2 + \mathcal{N}(\cdot,\cdot,p,\dot{p},v_y)\right](y,s) \Big| \rmd y \rmd s
\\ \nonumber & &
+ \int_0^t \int_\mathbb{R} \Big| [e(x,t-s+1) - e(x,t+1)] \psi(y) \left[ v_y^2 + \mathcal{N}(\cdot,\cdot,p,\dot{p},v_y)\right](y,s) \Big| \rmd y \rmd s.
\end{eqnarray}
In the remainder of this section, we will estimate the right-hand side of (\ref{phi-estimate}) term by term.

First, we note that (\ref{e:tildeg}) implies that there is a constant $C_1=C_1(M)$ so that $\tilde{\mathcal{G}}(x,0,t+1)\leq C_1\theta_1(x,t)$. Using (\ref{p-estimate}) and the fact that $|p(0)|\leq C_1\epsilon$, we therefore obtain
\[
(|p(0)| + |p(t)|) \tilde{\mathcal{G}}(x,0,t+1) \leq C_1 (\epsilon + h_2(t)) \theta_1(x,t),
\]
which is the desired estimate for the first term on the right-hand side of (\ref{phi-estimate}).

Next, we consider the integral term in (\ref{phi-estimate}) that involves the initial data $v_0$. Using (\ref{e:tildeg}) together with our assumption (\ref{e:phi0}) on $\phi_0$, and hence on $v_0$, we see that
\begin{equation}\label{est-initial}
\int_\mathbb{R} |\tilde{\mathcal{G}}(x,y,t) v_0(y)| \rmd y \leq C \epsilon \int_{\mathbb{R}} t^{-1/2} \left(\rme^{-\frac{(x-y+ct)^2}{4t}} + \rme^{-\frac{(x-y-ct)^2}{4t}} \right) \rme^{-\frac{y^2}{M}} \rmd y,
\end{equation}
which is clearly bounded by $\epsilon C_1\theta_1(x,t)$ for $t\geq1$ upon using
\[
\rme^{-\frac{(x-y\pm ct)^2}{4t}} \rme^{-\frac{y^2}{M}} \leq C_1 \rme^{-\frac{(x\pm ct)^2}{Mt}} \rme^{-\frac{y^2}{2M}}.
\]
For $t\leq1$, we can use the estimates
\begin{eqnarray*}
\rme^{-\frac{(x-y\pm ct)^2}{8t}} \rme^{-\frac{y^2}{M}} \;\leq\; 2\rme^{-\frac{(x-y)^2}{8t}}\rme^{-\frac{y^2}{M}} & \leq & C_1 \rme^{-\frac{x^2}{2M}} \\
\int_{\mathbb{R}} t^{-1/2} \left(\rme^{-\frac{(x-y+ct)^2}{8t}} + \rme^{-\frac{(x-y-ct)^2}{8t}} \right) \rmd y & \leq & C_1
\end{eqnarray*}
to conclude that the integral in (\ref{est-initial}) is again bounded by $\epsilon C_1 \theta_1(x,t)$.

We now consider the remaining two integrals in (\ref{phi-estimate}). Note that the definition of $h_1$ gives
\[
|v_y(y,s)| \leq \theta_2(y,s) h_1(s).
\]
Using this fact together with the estimates (\ref{p-estimate}) and (\ref{E:p_est}) and the bound (\ref{bound-N}), we obtain
\begin{eqnarray*}
|\mathcal{N}(y,s,p(s),\dot{p}(s),v_y(y,s))| & \leq & C_1 (1+s)^{-1/2} \left(\rme^{-\frac{(y + cs)^2}{8(s+1)}} + \rme^{-\frac{(y - cs)^2}{8(s+1)}} \right) |v_y(y,s)| |p(s)| \\ & &
+ C_1  \left(\rme^{-\frac{(y + cs)^2}{8(s+1)}} + \rme^{-\frac{(y - cs)^2}{8(s+1)}} \right)  |p(s)\dot{p}(s)|
\\ & \leq &
C_1 (\epsilon + h_2(t)) h_1(t) \theta_2(y,s) (1+s)^{-1/2} \left(\rme^{-\frac{(y + cs)^2}{8(s+1)}} + \rme^{-\frac{(y - cs)^2}{8(s+1)}} \right) \\ &&
+ C_1 (\epsilon + h(t)^2) \rme^{-\frac{c^2 s}{M}} \left(\rme^{-\frac{(y + cs)^2}{8(s+1)}} + \rme^{-\frac{(y - cs)^2}{8(s+1)}} \right) \\ & \leq &
C_1 (\epsilon + h(t)^2) \left[ (1+s)^{\gamma-1/2} \theta_1(y,s) \theta_2(y,s) + (1+s)^{\gamma} \theta_1(y,s) \rme^{-\frac{c^2 s}{M}} \right].
\end{eqnarray*}
In \S\ref{S:lemmas}, we will prove the following result for our spatio-temporal template functions.

\begin{Lemma}\label{lem-est-tG}
For each sufficiently large $M$, there is a constant $C_1$ so that
\begin{eqnarray*}
\int_0^t \int_\mathbb{R} |\tilde{\mathcal{G}}(x,y,t-s)| \left[ \theta_2^2 + (1+s)^{\gamma-1/2}\theta_1\theta_2+(1+s)^{\gamma} \theta_1\rme^{-c^2 s/M}\right](y,s) \rmd y \rmd s & \leq & C_1 \theta_1(x,t) \\
\int_0^t \int_\mathbb{R} |\tilde{\mathcal{G}}_x(x,y,t-s)| \left[\theta_2^2 + (1+s)^{\gamma-1/2}\theta_1\theta_2+(1+s)^{\gamma} \theta_1\rme^{-c^2 s/M}\right](y,s) \rmd y \rmd s & \leq & C_1\theta_2(x,t).
\end{eqnarray*}
\end{Lemma}

Using this lemma and the above estimates for $v_y^2+\mathcal{N}$, we obtain the desired estimate
\[
\int_0^t \int_\mathbb{R} \tilde{\mathcal{G}}(x,y,t-s)\left[ v_y^2 +\mathcal{N}(\cdot,\cdot,p,\dot{p},v_y)\right](y,s) \rmd y \rmd s \leq C_1 (\epsilon + h(t)^2) \theta_1(x,t).
\]
Finally, we have the following lemma, whose proof is again given in \S\ref{S:lemmas}, which provides the desired estimate for the last integral in (\ref{phi-estimate}).

\begin{Lemma}\label{lem-est-ediff}
For each sufficiently large $M$, there is a constant $C_1$ so that
\begin{eqnarray*}
\int_0^t \int_\mathbb{R} |e(x,t-s+1) - e(x,t+1)|\psi(y) \left[v_y^2 + \mathcal{N}(\cdot,\cdot,p,\dot{p},v_y) \right](y,s) \rmd y \rmd s & \leq & C_1 (\epsilon + h(t)^2)\theta_1(x,t) \\
\int_0^t \int_\mathbb{R} |e_x(x,t-s+1) - e_x(x,t+1)|\psi(y)\left[v_y^2 +\mathcal{N}(\cdot,\cdot,p,\dot{p},v_y)\right](y,s) \rmd y \rmd s & \leq & C_1(\epsilon + h(t)^2)\theta_2(x,t). 
\end{eqnarray*}
\end{Lemma}

In summary, combining the estimates obtained above, we have established the claimed estimate (\ref{key-est-phi}), and it remains to derive the estimate (\ref{key-est-phix}) to complete the proof of Proposition~\ref{P:claim-zeta}. Taking the $x$-derivative of equation (\ref{E:v_p}), we see that
\begin{eqnarray}\label{phix-estimate}
|v_x(x,t)| & \leq & \left(|p(0)| + |p(t)| \right) \tilde{\mathcal{G}}_x(x,0,t) + \int_\mathbb{R}\tilde{\mathcal{G}}_x(x,y,t) v_0(y)\rmd y \\ \nonumber &&
+ \int_0^t \int_\mathbb{R} \tilde{\mathcal{G}}_x(x,y,t-s) \left[ v_y^2 + \mathcal{N}(\cdot,\cdot,p,\dot{p},v_y) \right](y,s) \rmd y \rmd s
\\ \nonumber & &
+ \int_0^t \int_\mathbb{R} [e_x(x,t-s+1) - e_x(x,t+1)]\psi(y) \left[ v_y^2 + \mathcal{N}(\cdot,\cdot,p,\dot{p},v_y) \right](y,s) \rmd y \rmd s.
\end{eqnarray}
Applying the second estimate in Lemmas~\ref{lem-est-tG} and~\ref{lem-est-ediff} to (\ref{phix-estimate}) and using that $\tilde{\mathcal{G}}_x(x,0,t)\leq C\theta_2(x,t)$, we immediately obtain (\ref{key-est-phix}).

%%%%%%%%%%%%%%%%%%%%%%%%%%%%%%%%%%%%%%%%%%%%%%%%%

\subsection{Proofs of Lemmas~\ref{lem-est-tG} and~\ref{lem-est-ediff}}\label{S:lemmas}

It remains to prove the lemmas that we used in the preceding section.

\begin{proof}[\textbf{Proof of Lemma~\ref{lem-est-tG}}]
We need to show that for each large $M$ there is a constant $C_1$ so that
\[
\int_0^t \int_\mathbb{R} |\tilde{\mathcal{G}}(x,y,t-s)| \left[\theta_2^2 + (1+s)^{\gamma-1/2} \theta_1\theta_2 + (1+s)^{\gamma} \theta_1 \rme^{-c^2 s/M} \right](y,s) \rmd y \rmd s \leq C_1 \theta_1(x,t)
\]
for all $t\geq0$. First, we note that there are constants $C_1,\tilde{C}_1>0$ such that
\[
\tilde{C}_1 \rme^{-y^2/M} \leq |\theta_1(y,s)| + |\theta_2(y,s)| \leq C_1 \rme^{-y^2/M}
\]
for all $0\leq s\leq t\leq1$. Thus, for some constant $C_1$ that may change from line to line, we have
\begin{eqnarray*}
\lefteqn{\int_0^t \int_\mathbb{R} |\tilde{\mathcal{G}}(x,y,t-s)|\left[\theta_2^2 + (1+s)^{\gamma-1/2}\theta_1\theta_2+ (1+s)^{\gamma} \theta_1\rme^{-\frac{c^2 s}{M}}\right](y,s) \rmd y \rmd s}
\\ & \leq &
C_1 \int_0^t \int_\mathbb{R} (t-s)^{-1/2} \rme^{-\frac{(x-y)^2}{4(t-s)}} \rme^{-\frac{y^2}{M}} \rmd y\rmd s
\\ & \leq &
C_1 \int_0^t \left[\int_{\{|y|\ge 2|x|\}} (t-s)^{-1/2} \rme^{-\frac{(x-y)^2}{8(t-s)}} \rme^{-\frac{x^2}{8(t-s)}} \rmd y+ \int_{\{|y|\le 2|x|\}} (t-s)^{-1/2} \rme^{-\frac{(x-y)^2}{4(t-s)}} \rme^{-\frac{4x^2}{M}} \rmd y\right]\rmd s
\\ & \leq &
C_1 \int_0^t \Big[ \rme^{-\frac{x^2}{8(t-s)}}+ \rme^{-\frac{4x^2}{M}} \Big]\rmd s
\\ & \leq &
C_1 \rme^{-\frac{4x^2}{M}}
\\ & \leq &
\frac{C_1}{\tilde{C}_1} \theta_1(x,t)
\end{eqnarray*}
for all $0\leq t\leq1$. An analogous computation can be carried out for the $x$-derivative since $\int_0^t (t-s)^{-1/2}\rmd s$ is bounded uniformly in $0\leq t\leq1$.

Thus, it remains to estimate the expression
\[
\theta_1(x,t)^{-1} \int_0^t \int_\mathbb{R} |\tilde{\mathcal{G}}(x,y,t-s)| \left[ \theta_2^2 + (1+s)^{\gamma-1/2}\theta_1\theta_2 + (1+s)^{\gamma} \theta_1\rme^{-c^2 s/M}\right](y,s) \rmd y \rmd s
\]
for $t\geq1$. Combining only the exponentials in this expression, we obtain terms that can be bounded by
\begin{equation}\label{exp-form}
\exp\left(\frac{(x+\alpha_3 ct)^2}{M(1+t)} - \frac{(x-y+\alpha_1 c(t-s))^2}{4(t-s)} - \frac{(y+\alpha_2 cs)^2}{M(1+s)} \right)
\end{equation}
with $\alpha_j=\pm c$. To estimate this  expression, we proceed as in \cite[Proof of Lemma~7]{HowardZumbrun06} and complete the square of the last two exponents in (\ref{exp-form}). Written in a slightly more general form, we obtain
\begin{eqnarray*}
\lefteqn{\frac{(x-y-\alpha_1(t-s))^2}{M_1(t-s)} + \frac{(y-\alpha_2 s)^2}{M_2(1+s)} \;=\;
\frac{(x-\alpha_1(t-s)-\alpha_2 s)^2}{M_1(t-s)+M_2(1+s)} } \\ &&
+ \frac{M_1(t-s)+M_2(1+s)}{M_1M_2(1+s)(t-s)}\left( y - \frac{xM_2(1+s) - (\alpha_1M_2(1+s) + \alpha_2M_1s)(t-s)}{M_1(t-s)+M_2(1+s)}\right)^2
\end{eqnarray*}
and conclude that the exponent in (\ref{exp-form}) is of the form
\begin{eqnarray}\label{exp-form-2}
\lefteqn{ \frac{(x+\alpha_3t)^2}{M(1+t)} - \frac{(x-\alpha_1(t-s)-\alpha_2s)^2}{4(t-s)+M(1+s)} } \\ \nonumber &&
- \frac{4(t-s)+M(1+s)}{4M(1+s)(t-s)}\left( y - \frac{xM(1+s) - (\alpha_1M(1+s) + 4\alpha_2s)(t-s)}{4(t-s)+M(1+s)}\right)^2,
\end{eqnarray}
with $\alpha_j=\pm c$. Using that the maximum of the quadratic polynomial $\alpha x^2+\beta x+\gamma$ is $-\beta^2/(4\alpha)+\gamma$, it is easy to see that the sum of the first two terms in (\ref{exp-form-2}), which involve only $x$ and not $y$, is less than or equal to zero. Omitting this term, we therefore obtain the estimate
\begin{eqnarray}\label{est-on-exp}
\lefteqn{ \exp\left( \frac{(x\pm ct)^2}{M(1+t)} - \frac{(x-y\delta_1c(t-s))^2}{4(t-s)} - \frac{(y-\delta_2cs)^2}{M(1+s)} \right) } \\ \nonumber & \leq &
\exp \left( - \frac{4(t-s)+Ms}{4M(1+s)(t-s)}\left( y - \frac{xM(1+s)+c(\delta_1M(1+s)+ 4\delta_2s)(t-s)}{4(t-s)+M(1+s)}\right)^2 \right)
\end{eqnarray}
for $\delta_j=\pm1$. Using this result, we can now estimate the integral (\ref{exp-form}) term by term using the key assumption that $0<\gamma<\frac12$. The term involving $\theta_2^2$ can be estimated as follows using (\ref{est-on-exp}):
\begin{eqnarray*}
\lefteqn{ \theta_1(x,t)^{-1} \int_0^t \int_\mathbb{R}   |\tilde{\mathcal{G}}(x,y,t-s)|\theta^2_2(y,s) \rmd y \rmd s } \\ & \leq &
C_1 (1+t)^{\gamma}\int_0^t \frac{1}{\sqrt{t-s}(1+s)^{1+2\gamma}} \\ && \times 
\int_\R \exp\left(- \frac{4(t-s)+M(1+s)}{4M(1+s)(t-s)}\left( y - \frac{[xM(1+s) \pm c(M(1+s) + 4s)(t-s)]}{4(t-s)+M(1+s)}\right)^2 \right) \rmd y\rmd s \\ & \leq &
C_1(1+t)^{\gamma}\int_0^t \frac{1}{\sqrt{t-s}(1+s)^{1+2\gamma}}  \sqrt{\frac{4M(1+s)(t-s)}{4(t-s)+M(1+s)}} \rmd s \\ & \leq &
C_1(1+t)^{\gamma}\int_0^{t/2} \frac{1}{(1+s)^{1/2+2\gamma}}\frac{1}{(1+t)^{1/2}} \rmd s + C_1(1+t)^{\gamma}\int_{t/2}^t \frac{1}{(1+s)^{1+2\gamma}}\rmd s \\ & \leq &
C_1(1+t)^{\gamma-1/2}+ C_1(1+t)^{-\gamma},
 \end{eqnarray*}
which is clearly bounded since $\gamma<\frac12$. Similarly, we have
\begin{eqnarray*}
\lefteqn{\theta_1(x,t)^{-1} \int_0^t \int_\mathbb{R}   |\tilde{\mathcal{G}}(x,y,t-s)|(1+s)^{\gamma-1/2}\theta_1\theta_2 (y,s) \rmd y \rmd s } \\ & \leq &
C_1(1+t)^{\gamma}\int_0^t \frac{1}{\sqrt{t-s}(1+s)^{\gamma+1}}  \sqrt{\frac{4M(1+s)(t-s)}{4(t-s)+M(1+s)}} \rmd s \\ & \leq &
C_1(1+t)^{\gamma}\int_0^{t/2} \frac{1}{(1+s)^{\gamma+1/2}}\frac{1}{(1+t)^{1/2}} \rmd s + C_1(1+t)^{\gamma}\int_{t/2}^t \frac{1}{(1+t)^{\gamma+1/2}} \frac1{(1+s)^{1/2}}\rmd s
\\ & \leq &
C_1(1+t)^{\gamma-1/2} + C_1,
\end{eqnarray*}
which is again bounded due to $\gamma<\frac12$. Finally, we estimate
\begin{eqnarray*}
\lefteqn{\theta_1(x,t)^{-1} \int_0^t \int_\mathbb{R} |\tilde{\mathcal{G}}(x,y,t-s)| (1+s)^{\gamma}  \rme^{-\frac{c^2 s}{M}}\theta_1(y,s) \rmd y \rmd s} \\ \nonumber & \leq &
C_1(1+t)^{\gamma}\int_0^t \frac{\rme^{-\frac{c^2 s}{M}}}{\sqrt{t-s}}  \sqrt{\frac{4M(1+s)(t-s)}{4(t-s)+M(1+s)}} \rmd s \\ \nonumber & \leq &
C_1(1+t)^{\gamma}\int_0^{t/2}\rme^{-\frac{c^2 s}{M}}\frac{1}{(1+t)^{1/2}} \rmd s + C_1(1+t)^{\gamma} \rme^{-\frac{c^2 t}{2M}} \int_{t/2}^t\rmd s  \\ \nonumber & \leq &
C_1(1+t)^{\gamma-1/2} + C_1(1+t)^{\gamma+1} \rme^{-\frac{c^2 t}{2M}},
\end{eqnarray*}
which is bounded, again due to $\gamma<\frac12$.

It remains to verify the second inequality in Lemma~\ref{lem-est-tG} which involves $\tilde{\mathcal{G}}_x$. We shall check only the term involving $(1+s)^{\gamma-1/2}\theta_1\theta_2$ as the other cases are similar and, in fact, easier. We have shown above that the resulting integrals are bounded for $0\leq t\leq1$ and therefore focus on the case $t\geq1$. Using that $|\tilde{\mathcal{G}}_x|\leq Ct^{-1/2}|\tilde{\mathcal{G}}|$, which follows by inspection, and employing again (\ref{est-on-exp}), we obtain
\begin{eqnarray*}
\lefteqn{\theta_2(x,t)^{-1} \int_0^t \int_\mathbb{R}   |\tilde{\mathcal{G}}_x(x,y,t-s)|(1+s)^{\gamma-1/2}\theta_1\theta_2(y,s) \rmd y \rmd s} \\ & \leq &
C_1(1+t)^{\gamma+1/2}\int_0^t \frac{1}{(t-s)(1+s)^{\gamma+1}} \sqrt{\frac{4M(1+s)(t-s)}{4(t-s)+M(1+s)}} \rmd s \\ & \leq &
C_1(1+t)^{\gamma+1/2}\int_0^{t/2} \frac{1}{t^{1/2}(1+s)^{\gamma+1/2}}\frac{1}{(1+t)^{1/2}}\rmd s + C_1\int_{t/2}^t\frac{1}{(t-s)^{1/2}} \frac{1}{(1+t)^{1/2}}\rmd s
\\ & \leq &
C_1(1+t)^{\gamma}t^{-1/2} + C_1,
\end{eqnarray*}
which is bounded for $t\geq1$. This completes the proof of Lemma~\ref{lem-est-tG}.
\end{proof}

\begin{proof}[\textbf{Proof of Lemma~\ref{lem-est-ediff}}]
We need to show that
\[
\int_0^t \int_\mathbb{R} |e(x,t-s+1) - e(x,t+1)|\psi(y) [v_y^2  +\mathcal{N}(\cdot,\cdot,p,\dot{p},v_y)](y,s) \rmd y \rmd s \leq C_1 (\epsilon + h(t)^2)\theta_1(x,t).
\]
Intuitively, this integral should be small for the following reason. The difference $e(x,t-s)-e(x,t+1)$ converges to zero as long as $s$ is not too large, say on the interval $s\in[0,t/2]$. For $s\in[t/2,t]$, on the other hand, we will get exponential decay in $s$ from the localization of $\psi(y)$ in combination with the propagating Gaussians that appear in the nonlinearity and forcing terms. To make this precise, we write
\begin{equation}\label{e-difference}
e(x,t-s) - e(x,t+1) = \underbrace{e(x,t-s) - e(x,t-s+1)}_{\mathrm{term~I}} + \underbrace{e(x,t-s+1) - e(x,t+1)}_{\mathrm{term~II}}.
\end{equation}
We focus first on the term I and consider the cases $t\geq1$ and $0\leq t\leq1$ separately. First, let $t\geq1$. For $0\leq s\leq t-1$, we have $|e(x,t-s)-e(x,t-s+1)|\leq C\tilde{\mathcal{G}}(x,0,t-s)$, and we can estimate the resulting integral above in the same way as in the proof of Lemma~\ref{lem-est-tG}; we omit the details. For $t-1\leq s\leq t\leq1$, on the other hand, the definition of $e(x,t-s)$ yields
\begin{equation}\label{est-t-small}
|e(x,t-s) - e(x,t-s+1)| \leq C \int_{\frac{x^2}{(1+t-s)}}^{\frac{x^2}{(t-s)}} \rme^{-z^2} \rmd z \leq C \rme^{-x^2/2}.
\end{equation}
Using (\ref{est-phiN})-(\ref{est-phiN3}), namely
\begin{equation}\label{est-N13}
|\psi(y)\mathcal{N}(y,s,p,\dot{p},v_y)| \leq C\rme^{-\frac{c}{2}|y| - \frac{2c^2}{M}s} (\epsilon + h(t)^2)
\end{equation}
for all $s\geq0$, we obtain
\begin{eqnarray*}
\lefteqn{\int_{t-1}^t \int_\mathbb{R} [e(x,t-s) - e(x,t-s+1)]\psi(y)|\mathcal{N}(y,s,p,\dot{p},v_y)|(y,s)\rmd y \rmd s} \\
& \leq & C_1 (\epsilon + h(t)^2) \int_{t-1}^t \rme^{-\frac{x^2}{2}}\rme^{-\frac{2c^2 s}{M}}\rmd s
\;\leq\; C_1 (\epsilon + h(t)^2) \rme^{-\frac{x^2}{2}}\rme^{-\frac{2c^2 t}{M}},
\end{eqnarray*}
which is clearly bounded by $C_1\theta_1(x,t)$ since $\rme^{-c^2 t/M}\leq C_1(1+t)^{-\gamma}$ and
\[
\frac {(x+ ct)^2}{M(1+t)} \leq \frac{2x^2}{M} + \frac{4c^2}{M} + \frac{c^2t}{M}
\]
for arbitrary $M\geq4$. In summary, we have established the desired estimates for the term I in (\ref{e-difference}) for $t\geq1$. For $t\leq1$, the estimate (\ref{est-t-small}) remains true since $t-s$ is small, and proceeding as above yields
\[
\int_{0}^t \int_\mathbb{R} [e(x,t-s) - e(x,t-s+1)]\psi(y) |\mathcal{N}(y,s,p,\dot{p},v_y)| \rmd y \rmd s \leq C(\epsilon  + h(t)^2) \rme^{-\frac{x^2}{2}}\rme^{-\frac{2c^2 t}{M}},
\]
which is again bounded by $C_1\theta_1(x,t)$. 

It remains to discuss the term II which involves the difference $e(x,t-s+1)-e(x,t+1)$. We have
\begin{eqnarray*}
\lefteqn{|e(x,t-s+1) - e(x,t+1)|} \\ & = &
\left|\int_{t+1}^{t-s+1} e_\tau(x,\tau)\rmd \tau\right| \\ & \leq &
\int_{t-s+1}^{t+1} \left| \frac{c}{\sqrt{4\pi \tau}} \left( \rme^{-\frac{(x-c\tau)^2}{4\tau}} + \rme^{-\frac{(x+c\tau)^2}{4\tau}}\right) + \frac{1}{\tau\sqrt{4\pi}} \left( \frac{(x-c\tau)}{\sqrt{4\tau}}\rme^{-\frac{(x-c\tau)^2}{4\tau}} - \frac{(x+c\tau)}{\sqrt{4\tau}}\rme^{-\frac{(x+c\tau)^2}{4\tau}}\right)   \right|\rmd \tau
\\ & \leq &
C\int_{t-s+1}^{t+1} \left( \frac{1}{\sqrt{\tau}} +  \frac{1}{\tau}  \right)  \left( \rme^{-\frac{(x-c\tau)^2}{8\tau}} + \rme^{-\frac{(x+c\tau)^2}{8\tau}}\right) \rmd \tau,
\end{eqnarray*}
where we used in the last inequality that $z\rme^{-z^2}$ is uniformly bounded in $z$. We now use the preceding expression to estimate $\theta_1^{-1}(x,t)(e(x,t-s+1)-e(x,t+1))$ and focus first on the single exponential term
\[
\rme^{\frac{(x-ct)^2}{M(1+t)}}\rme^{-\frac{(x-c\tau)^2}{8\tau}}.
\]
Combining these exponentials and completing the square in $x$ in the resulting exponent, the latter becomes
\[
-\frac{[M(t-\tau+1) + (M-8)\tau]}{8M(t+1)\tau}\left[ x + \frac{c(8-M)\tau(t+1)}{M(t-\tau) + (M-8)\tau}\right]^2 + \frac{c^2(t-\tau+1)^2}{M(t-\tau+1) + (M-8)\tau}.
\]
Using that $\tau\leq t$ and picking $M\geq8$, we can neglect the exponent resulting from the first expression that involves in $x$ and conclude that
\[
\rme^{\frac{(x-ct)^2}{M(1+t)}}\rme^{-\frac{(x-c\tau)^2}{8\tau}} \leq C_1 \rme^{\frac{c^2(t-\tau)}{M}}.
\]
The remaining exponentials can be estimated similarly, and we obtain
\begin{eqnarray*}
\theta_1^{-1}(x,t) |e(x,t-s+1) - e(x,t+1)| & \leq &
C_1 (1+t)^{\gamma}\int_{t-s+1}^{t+1} \left( \frac{1}{\sqrt{\tau}} +  \frac{1}{\tau}  \right) \rme^{\frac{c^2(t-\tau)}{M}} \rmd \tau \\ \nonumber & \leq & 
C_1 (1+t)^{\gamma} (1+t-s)^{-1/2}\rme^{\frac{c^2s}{M}}.
\end{eqnarray*}
Using this inequality together with (\ref{est-N13}) finally gives
\begin{eqnarray*}
\lefteqn{\theta_1(x,t)^{-1}\int_0^t \int_\mathbb{R} [e(x,t-s+1) - e(x,t+1)]\psi(y)|\mathcal{N}(y,s,p,\dot{p},v_y)|(y,s) \rmd y \rmd s} \\ & \leq &
C_1(1+t)^{\gamma} (\epsilon + h(t)^2)\int_0^t (1+t-s)^{-1/2}\rme^{\frac{c^2s}{M}} \rme^{-\frac{2c^2 s}{M}}\rmd s \\ & \leq &
C_1(1+t)^{\gamma} (\epsilon + h(t)^2)\left[(1+t)^{-1/2}\int_0^{t/2}  \rme^{-\frac{c^2 s}{M}}\rmd s + \rme^{-\frac{c^2 t}{2M}}\int_{t/2}^t (1+t-s)^{-1/2}\rmd s\right] \\ & \leq &
C_1(\epsilon + h(t)^2).
\end{eqnarray*}
for $M$ sufficiently large, which proves the first estimate in Lemma~\ref{lem-est-ediff}.

It remains to prove the estimate
\begin{eqnarray*}
&& \int_0^t \int_\mathbb{R} |e_x(x,t-s+1) - e_x(x,t+1)|\psi(y)\left[v_y^2 +\mathcal{N}(\cdot,\cdot,p,\dot{p},v_y)\right](y,s) \rmd y \rmd s \\ && \qquad\qquad \leq C_1(\epsilon + h(t)^2)\theta_2(x,t)
\end{eqnarray*}
for the derivative in $x$. Since the derivative of $e(x,t-s+1)-e(x,t+1)$ with respect to $x$ generates an extra decay term $(1+t)^{-1/2}$, we have
\begin{eqnarray*}
\lefteqn{ \theta_2(x,t)^{-1}\int_0^t \int_\mathbb{R}  [e_x(x,t-s+1) - e_x(x,t+1)]\psi(y) |\mathcal{N}(y,s,p,\dot{p},v_y)|(y,s) \rmd y \rmd s } \\ & \leq &
C_1 (\epsilon ^2 + h(t)^2)(1+t)^{\gamma+1/2}\int_0^t (1+t-s)^{-1}\rme^{\frac{c^2s}{M}} \rme^{-\frac{2c^2 s}{M}} \rmd s \\ & \leq & C_1 (\epsilon ^2 + h(t)^2)(1+t)^{\gamma+1/2}\left[(1+t)^{-1}\int_0^{t/2}\rme^{-\frac{c^2 s}{M}} \rmd s + \rme^{-\frac{c^2 t}{2M}}\int_{t/2}^t(1+t-s)^{-1}\rmd s\right] \\ & \leq & C_1 (\epsilon ^2 + h(t)^2),
\end{eqnarray*}
which completes the proof of the lemma. 
\end{proof}

%%%%%%%%%%%%%%%%%%%%%%%%%%%%%%%%%%%%%%%%%%%%%%%%%%%%%%%%%%%%%%%%%%%%%%%%%%%%

\paragraph{Acknowledgments.}
Beck, Nguyen, Sandstede, and Zumbrun were supported partially by the NSF through grants DMS-1007450, DMS-1108821, DMS-0907904, and DMS-0300487, respectively.

%%%%%%%%%%%%%%%%%%%%%%%%%%%%%%%%%%%%%%%%%%%%%%%%%%%%%%%%%%%%%%%%%%%%%%%%%%%%

\bibliographystyle{sandstede}
%\bibliography{toy_model}

\end{document}